\documentclass[12pt]{amsart}

\usepackage{amsmath,amssymb,amsthm}
\usepackage{amscd,times,fullpage}
\newtheorem{thm}{Theorem}
\newtheorem{prop}{Proposition}
\newtheorem{lem}{Lemma}
\newtheorem{cor}{Corollary}

\theoremstyle{remark}
\newtheorem{rem}{Remark}

\theoremstyle{definition}
\newtheorem{defn}{Definition}
\newtheorem{ex}{Example}

\newcommand{\C}{\mathbb{ C}}

\newcommand{\N}{\mathbb{ N}}

\newcommand{\Q}{\mathbb{ Q}}
\newcommand{\R}{\mathbb{ R}}

\newcommand{\rk}{\operatorname{rk}}

\newcommand{\Sym}{\operatorname{Sym}}
\newcommand{\Sh}{\operatorname{Sh_\lambda}}
\newcommand{\Z}{\mathbb{ Z}}

\newcommand{\sg}{\mathrm{sign}}
\newcommand{\rw}{\mathrm{w}}

\DeclareMathOperator{\dummygg}{\mathfrak{g}}
\renewcommand{\gg}{\dummygg}
\DeclareMathOperator{\II}{\mathfrak{I}}
\DeclareMathOperator{\hh}{\mathfrak{h}}
\DeclareMathOperator{\TT}{\mathfrak{t}}

\DeclareMathOperator{\LL}{\mathfrak{L}}
\DeclareMathOperator{\ad}{ad}
\DeclareMathOperator{\Ad}{Ad}

\title{Complex cobordism classes of homogeneous spaces}
\author{Victor M.~Buchstaber}
\address{Steklov Mathematical Institute, Russian Academy
of Sciences, Gubkina Street 8, 119991 Moscow, Russia}
\email{buchstab@mi.ras.ru}
\address{School of Mathematics, University of Manchester, Oxford Road,
Manchester M13~9PL, England}
\email{Victor.Buchstaber@manchester.ac.uk}
\author{Svjetlana Terzi\'c}
\address{Faculty of Science, University of Montenegro,
Cetinjski put bb, 81000 Podgorica, Montenegro}
\email{sterzic@cg.ac.yu}
\date{\today ; MSC 2000: primary 57R77, 22F30; secondary 57R20, 14M15 }

\begin{document}

\maketitle

\begin{abstract}
We consider compact homogeneous spaces $G/H$ of positive Euler
characteristic endowed with an invariant almost complex structure
$J$ and the canonical action $\theta$  of the maximal torus $T ^{k}$
on $G/H$. We obtain explicit formula for the cobordism class of such
manifold through the weights of the action  $\theta$ at the identity
fixed point $eH$ by an action of the quotient group $W_G/W_H$ of the
Weyl groups for $G$ and $H$. In this way we show that the cobordism
class for such manifolds can be computed explicitly without
information on their cohomology. We also show that formula for
cobordism class  provides an explicit way for computing the
classical Chern numbers for $(G/H, J)$.  As a consequence we obtain
that the Chern numbers for $(G/H, J)$ can be computed without
information on cohomology for $G/H$. As an application we provide an
explicit formula for cobordism classes and characteristic numbers of
the flag manifolds $U(n)/T^n$, Grassmann manifolds
$G_{n,k}=U(n)/(U(k)\times U(n-k))$ and some particular interesting
examples.

\end{abstract}

\section{introduction}
In this paper we consider the problem of description of  complex
cobordism classes of homogeneous spaces $G/H$ endowed with an
invariant almost complex structure, where $G$ is compact connected
Lie group and $H$ is its closed connected subgroup of maximal rank.
These spaces are classical manifolds and have a very reach geometric
structure from the different points of view and our interest in
these manifolds is related to the well known  problem in cobordism
theory to find the representatives  in cobordism classes that have
reach geometric structure. Our interest in research of the
homogeneous spaces $G/H$ with positive Euler characteristic,  is
also stimulated  by well known relations between the cohomology
rings of these spaces and the deep problems in the theory of
representations and combinatorics (see, for example~\cite{F}). These
problems are formulated in terms of different additive basis in
cohomology rings for $G/H$ and multiplicative rules related to that
basis. We hope the   research of the cobordisms of the spaces $G/H$
to give the new relations and bring the new results in that
direction.

We use the approach based on the notion of Chern-Dold character
originally introduced in~\cite{Buchstaber} and the notion of
universal toric genus introduced in~\cite{BR} and described in
details in~\cite{IMRN}.The universal toric genus can be constructed
for any even dimensional manifold $M^{2n}$ with a given torus action
and stable complex structure which is equivariant under given torus
action. Moreover, if the set of isolated fixed points for this
action is finite that the universal toric genus can be localized,
which means that it can be explicitly written through the weights
and the signs at the fixed points for the representations that gives
arise from the given torus action.

The construction of the toric genus is reduced to the computation of
Gysin homomorphism of $1$ in complex cobordisms for fibration whose
fiber is $M^{2n}$ and the base is classifying space of the torus.
The problem of the localization of Gysin homomorphism is very known
and it was studied by many authors, starting with 60-es of the last
century. In~\cite{BR} and~\cite{IMRN} is obtained explicit answer
for this problem in the terms of the torus action on tangent bundle
for $M^{2n}$. The history of this problem is presented also in these
papers.

If consider homogeneous space $G/H$ with $\rk G = \rk H = k$, then
we have on it the canonical action $\theta$ of the maximal torus
$T^{k}$ for $H$ and $G$, and any invariant almost complex structure
on $G/H$ is compatible with this action. Besides that, all fixed
points for the action $\theta$ are isolated, so one can apply
localization formula to compute universal toric genus for this
action and any invariant almost complex structure. Since, in this
case, we consider almost complex structures, all fixed points in the
localization formula  are going to have sign $+1$. We prove that the
weights for the action $\theta$ at different fixed points can be
obtained  by an action of the Weyl group $W_{G}$ up to an action of
the Weyl group $W_{H}$ on the weights for $\theta$ at identity fixed
point. On this way we get an  explicit formula for the cobordism
classes of such spaces in terms of the weights at the fixed point
$eH$. This formula also shows that the cobordism class for $G/H$
related to an invariant almost complex structure can be computed
without information about cohomology for $G/H$.

 We obtain also the explicit formulas, in terms of the weights at
identity fixed point, for the cohomology characteristic numbers for
homogeneous spaces of positive Euler characteristic endowed with an
invariant almost complex structure. We use further that the
cohomology characteristic numbers $s_\omega, \; \omega=(i_1, \ldots,
i_n)$, and classical Chern numbers $c^\omega=c_1^{i_1}\cdots
c_n^{i_n}$ are related by some standard relations from the theory of
symmetric polynomials. This fact together with the obtained formulas
for the characteristic numbers $s_\omega(\tau(G/H))$ proves that the
classical Chern numbers $c^\omega(\tau(G/H))$ for the homogeneous
spaces under consideration can be computed without information on
their cohomology. It also gives an explicit way for the computation
of the classical Chern numbers.

We provide an application of our results by obtaining explicit
formula for the cobordism class and top cohomology characteristic
number of the flag manifolds $U(n)/T^n$ and Grassmann manifolds
$G_{n,k}=U(n)/(U(k)\times U(n-k))$ related to the standard complex
structures. We want to emphasize that, our method when applying to
the flag manifolds and Grassmann manifolds gives the description of their
cobordism classes and characteristic numbers using the technique of {\it divided
difference operators}. Our method also makes possible to compare
cobordism classes that correspond to the different invariant almost
complex structures on the same homogeneous space. We illustrate that
on  the space  $U(4)/(U(1)\times U(1)\times U(2))$, which is firstly
given in~\cite{BH} as an example of homogeneous space  that admits
two different invariant complex structure.

This paper comes out from the first part of our work where we mainly
considered invariant almost complex structures on homogeneous spaces
of positive Euler characteristic. It has continuation which is going
to deal with the same questions, but related to the  stable complex
structures equivariant under given torus action on homogeneous
spaces of positive Euler characteristic.

The authors are grateful to A.~Baker and A.~Gaifullin for useful
comments.

\section{Universal toric genus}
We will recall the results from~\cite{BR},~\cite{BPR} and~\cite{IMRN}.

\subsection{General setting.}In general setting one considers
$2n$-dimensional manifold $M^{2n}$ with a given smooth action
$\theta$ of the torus $T^{k}$. We say that $(M^{2n}, \theta,
c_{\tau})$ is {\it tangentially stable complex} if it admits
$\theta$-equivariant stable complex structure $c_{\tau}$. This means
that there exist $l\in \N$ and complex vector bundle $\xi$ such that
\begin{equation}\label{scdp}
c_{\tau}\colon \tau (M^{2n})\oplus \R ^{2(l-n)}\longrightarrow \xi
\end{equation}
is real isomorphism and the composition
\begin{equation}\label{scd}
  r(t)\colon\xi\stackrel{c_\tau^{-1}}{\longrightarrow}\tau(M^{2n})\oplus
  \mathbb{R}^{2(l-n)}\stackrel{d\theta(t)\oplus I}{\longrightarrow}
  \tau(M^{2n})\oplus\mathbb{R}^{2(l-n)}\stackrel{c_\tau}{\longrightarrow}\xi \
\end{equation}
is a complex transformation for any $t\in T^{k}$.

If there exists $\xi$ such that $c_{\tau}\colon \tau
(M^{2n})\longrightarrow \xi$ is an isomorphism, i.~e.~ $l=n$, then
$(M^{2n},\theta ,c_{\tau})$ is called {\it almost complex}
$T^k$-manifold.

Denote by $\Omega _{U}^{*} [[u_1,\ldots, u_{k}]]$ an algebra of
formal power series over $\Omega _{U}^{*} = U^{*}(pt)$. It is well
known~\cite{Novikov} that $U^{*}(pt) = \Omega _{U}^{*} = \Z
[y_{1},\ldots ,y_{n},\ldots ]$, where $\dim y_{n}=-2n$. Moreover, as
the generators for $\Omega _{U}^{*}$ over the rationales, or in
other words for $\Omega _{U}^{*}\otimes \Q$, can be taken the family
of cobordism classes $[\C P^{n}]$ of the complex projective spaces.

When given a $\theta$-equivariant stable complex structure
$c_{\tau}$ on $M^{2n}$,  we can always choose $\theta$-equivariant
embedding $i\colon M^{2n}\to \R ^{2(n+m)}$, where $m>n$, such that
$c_{\tau}$ determines, up to natural equivalence, a
$\theta$-equivariant complex structure $c_{\nu}$ on the normal
bundle $\nu(i)$ of $i$. Therefore, one can define the universal
toric genus for $(M^{2n},\theta, c_{\tau})$ in complex cobordisms,
see~~\cite{BR},~\cite{IMRN}.

We want to note that, in the case when $c_{\tau}$ is almost complex
structure, an universal toric genus for $(M^{2n},\theta ,c_{\tau})$ is
completely defined in terms of the action $\theta$ on tangent bundle
$\tau (M^{2n})$.

The universal toric genus for $(M^{2n}, \theta, c_{\tau})$ could be
looked at as an element in algebra $\Omega _{U}^{*} [[u_1,\ldots,
u_{k}]]$. It is defined with
\begin{equation}
\Phi (M^{2n}, \theta, c_{\tau}) =  [M^{2n}] + \sum _{|\omega|>0}
[G_{\omega}(M^{2n})]u^{\omega} \; ,
\end{equation}
where $\omega = (i_1,\ldots ,i_k)$ and $u^{\omega} =
u_{1}^{i_1}\cdots u_{k}^{i_k}$.

Here by $[M^{2n}]$ is denoted the complex cobordism class of the
manifold $M^{2n}$  with stable complex structure $c_{\tau}$, by
$G_{\omega}(M^{2n})$ is denoted the stable complex manifold obtained
as the total space of the fibration $G_{\omega} \to B_{\omega}$ with
fiber $M$. The base $B_{\omega} = \prod _{j=1}^{k}B_{j}^{i_j}$,
where $B_{j}^{i_j}$ is Bott tower, i.~e.~$i_j$-fold iterated
two-sphere bundle over $B_{0}=pt$. The base $B_{\omega}$ satisfies
$[B_{\omega}]=0$, $|\omega|>0$, where $|\omega| =
\sum\limits_{j=1}^{k} i_{j}$.

\subsection{The action  with isolated fixed points.}
Assume that the action of $T^{k}$ on $M^{2n}$ has isolated fixed points.
We first introduce, following~\cite{BR}, the general notion of the
{\it sign} at isolated fixed point. Let $p$ be an isolated fixed point.
The representation $r_{p}\colon T^k\to GL(l,\C)$ associated
to~\eqref{scd} produces the decomposition  of the fiber $\xi
_{p}\cong \C ^{l}$ as $\xi _{p}\cong \C ^{l-n}\oplus \C ^{n}$. In
this decomposition $r_{p}$ acts trivially on $\C ^{l-n}$ and without
trivial summands on $\C ^{n}$. From the other hand the isomorphism
$c_{\tau ,p}$ from~\eqref{scdp} defines an orientation in the
tangent space $\tau _{p}(M)$. This together leads to the following
definition.

\begin{defn}
 The  $\sg (p)$ at isolated fixed point $p$ is
$+1$ if the map
\[
\tau_p(M^{2n})\stackrel{I\oplus\hspace{.1ex}0}{\longrightarrow}
\tau_p(M^{2n})\oplus \mathbb{R}^{2(l-n)} \stackrel{c_{\tau,p}}
{\longrightarrow}\xi_p \cong\mathbb{C}^n \oplus\mathbb{C}^{l-n}
\stackrel{\pi}{\longrightarrow} \mathbb{C}^n \; ,
\]
preserves orientation. Otherwise, $\sg (p)$ is $-1$.
\end{defn}
\begin{rem}\label{cmx}
Note that for an almost complex $T^{k}$-manifold $M^{2n}$, it
directly follows from the definition that  $\sg (p)= +1$ for any
isolated fixed point.
\end{rem}
If an action $\theta$ of $T^{k}$ on $M^{2n}$ has only isolated fixed
points, then it is proved that toric genus for $M^{2n}$  can be
completely described using just local data at the fixed
points,~\cite{BR},~\cite{IMRN}.

Namely, let $p$ again be an isolated fixed point. Then the non
trivial summand of $r_p$ from~(\ref{scd}) gives rise to the
tangential representation of $T^{k}$ in $GL(n, \C)$. This
representation decomposes into $n$ non-trivial one-dimensional
representations of $r_{p,1}\oplus \ldots \oplus r_{p,n}$ of $T^{k}$.
Each of the representations $r_{p,j}$ can be written as
\[
r_{p,j}(e^{2\pi i x_{1}}, \ldots ,e^{2\pi i x_{k}})v= e^{2\pi i
\langle\Lambda_{j}(p), {\bf x}\rangle }v \; ,
\]
for some $\Lambda_{j}(p) = (\Lambda^{1}_{j}(p),\ldots
,\Lambda^{k}_{j} (p)) \in \Z ^{k}$, where ${\bf x}=(x_{1},\ldots
,x_{k})\in \R ^{k}$ and $\langle \Lambda_{j}(p), {\bf x}\rangle =
\sum\limits _{l=1}^{k}\Lambda^{l}_{j}(p)x_{l}$. The sequence $\{
\Lambda_{1}(p),\ldots ,\Lambda_{n}(p)\}$ is called {\it the weight
vector} for representation $r_{p}$ in the fixed point $p$.

Let $F(u, v)= u + v + \sum \alpha _{ij}u^{i}v^{j}$ be {\it the
formal group for complex cobordism} \cite{Novikov-67}. The
corresponding power system $\{ [w](u)\in \Omega ^{*}[[u]] : w\in \Z
\}$ is uniquely defined with $[0](u)=0$ and $[w](u) = F(u,
[w-1])(u)$, for $w\in \Z$. For ${\bf w}=(w_1,\ldots ,w_{k})\in
{\Z}^{k}$ and ${\bf u}=(u_1,\ldots, u_k)$ one defines $\bf{[w](u)}$
inductively with ${\bf [w](u)} =[w](u)$ for $k=1$ and
\[
{\bf[w](u)} = F_{q=1}^{k}[w_{q}](u_{q})=F(F_{q=1}^{k-1}[w_{q}](u_{q}),
[w_{k}](u_{k})) \; ,
\]
for $k\geqslant 2$. Then for toric genus of the action $\theta$
with isolated fixed points the following localization formula holds,
which is first formulated in~\cite{BR} and proved in details
in~\cite{IMRN}.

\begin{thm}\label{EGFP}
If the action $\theta$ has a finite set $P$  of isolated fixed
points then
\begin{equation}
\Phi (M^{2n}, \theta, c_{\tau}) = \sum _{p\in P} \sg (p) \prod
_{j=1}^{n}\frac{1}{[\Lambda_{j}(p)]({\bf u})}
\end{equation}
and it is equal to $[M^{2n}] + \LL ({\bf u})$, where $\LL ({\bf u})
\in \Omega _{U}^{*} [[u_{1},\ldots ,u_{k}]]$ and $\LL ({\bf 0}) =
0$.
\end{thm}

\subsection{Chern-Dold character.} In review of the basic definitions
and results on Chern character we follow~\cite{Buchstaber}.

Let $U^{*}$ be the theory of unitary cobordisms.

\begin{defn}
The Chern-Dold character  for a topological space $X$ in the theory
of unitary cobordisms $U^{*}$ is a ring homomorphism
\begin{equation}
ch _{U} : U^{*}(X) \to H^{*}(X, \Omega _{U}^{*}\otimes \Q) \ .
\end{equation}
\end{defn}

Recall that the Chern-Dold character as a multiplicative
transformation of cohomology theories is uniquely defined by the
condition that for $X=(pt)$ it gives canonical inclusion $\Omega
_{U}^{*}\to \Omega _{U}^{*}\otimes \Q$.

The Chern-Dold character splits into composition

\begin{equation}\label{CD}
ch _{U} : U^{*}(X)\to H^{*}(X, \Omega _{U}^{*}(\Z ))\to H^{*}(X,
\Omega _{U}^{*}\otimes \Q ) \; .
\end{equation}

The ring $\Omega _{U}^{*}(\Z ))$ in \eqref{CD} is firstly described
in~\cite{Buchstaber}. It is a subring of $\Omega _{U}^{*}\otimes \Q$
generated by the elements from $\Omega _{U}^{-2n}\otimes \Q$ having
integers Chern numbers. It is equal to
\[
\Omega _{U}^{*}(\Z )=\Z [b_1,\ldots , b_n,\ldots ] \; ,
\]
where $b_n = \frac{[\C P^{n}]}{n+1}$.

Then the Chern character leaves  $[M^{2n}]$ invariant, i.~e.~
\begin{equation}
ch_{U}([M^{2n}])= [M^{2n}]\; ,
\end{equation}
and $ch_{U}$ is the homomorphism of $\Omega _{U}^*$-modules

It follows from the its description~\cite{Buchstaber} that the
Chern-Dold character $ch _{U} : U^{*}(X)\to H^{*}(X, \Omega
_{U}^{*}(\Z ))$ as a multiplicative transformation of the cohomology
theories is given by the series
\[
ch _{U}u = h(x)=\frac{x}{f(x)},\quad \mbox{where}\quad f(x)= 1+\sum
_{i=1}^{\infty} a_{i}x^{i} \quad \mbox{and}\quad a_{i}\in \Omega
_{U}^{-2i}(\Z) \; .
\]
Here $u=c_1^U(\eta) \in U^2(\mathbb{C}P^\infty), \; x=c_1^H(\eta)
\in H^2(\mathbb{C}P^\infty,\mathbb{Z})$ denote the first Chern classes
of the universal complex line bundle $\eta \to \mathbb{C}P^\infty$.

From the construction of Chern-Dold character follows also the
equality
\begin{equation}\label{charact}
ch_{U}([M^{2n}]) = [M^{2n}] = \sum _{\|\omega \|=n} s_{\omega}(\tau
(M^{2n}))a^{\omega}\; ,
\end{equation}
where $\omega = (i_{1},\ldots ,i_{n}), \; \| \omega \| = \sum
_{l=1}^{n} l\cdot i_{l}$ and $a ^{\omega} = a_{1}^{i_1}\cdots
a_{n}^{i_{n}}$. Here the numbers $s_{\omega}(\tau
(M^{2n})),{\|\omega \|=n}$ are {\it the cohomology characteristic
numbers} of $M^{2n}$.

If on $M^{2n}$ is given torus action $\theta$ of $T^{k}$ and stable
complex structure $c_{\tau}$ which is $\theta$-equivariant, then the
Chern character of its toric genus is

\begin{equation}\label{cheq}
ch_{U}(\Phi (M^{2n},\theta, c_{\tau} )) = [M^{2n}] + \sum _{\|
\omega \|
>0}[G_{\omega}(M^{2n})](ch_{U}{\bf u})^{\omega}\; ,
\end{equation}
where $ch_{U}{\bf u} = (ch_{U}u_1, \ldots, ch_{U}u_k)$ and
$ch_{U}u_i = \frac{x_i}{f(x_i)}$.

We have that $F(u, v)=g^{-1}\left( g(u)+g(v)\right)$, where $g(u)=u
+ \sum\limits_{n>0}\frac{[\C P^{n}]}{n+1}u^{n+1}$ (see
\cite{Novikov-67}) is {\it the logarithm of the formal group} $F(u,
v)$ and $g^{-1}(u)$ is the function inverse to the series $g(u)$.
Using that $ch_{U}g(u)=x$ (see \cite{Buchstaber}), we obtain
$ch_{U}F(u_1, u_2)=h(x_1+x_2)$ and therefore
\[ ch_{U}[\Lambda _{j}(p)](u) = \frac{\langle \Lambda _{j}(p),
{\bf x}\rangle}{f(\langle \Lambda _{j}(p), {\bf
x}\rangle)}\; .  \]

Applying these results to the theorem~\eqref{EGFP} we get

\begin{equation}\label{chs}
ch_{U}(\Phi(M^{2n},\theta, c_{\tau } ))= \sum _{p\in P}\sg (p)\prod
_{j=1}^{n}\frac{f(\langle \Lambda _{j}(p), {\bf x}\rangle)}{\langle
\Lambda _{j}(p),{\bf x}\rangle} \; .
\end{equation}

The formulas~\eqref{cheq} and~\eqref{chs} gives that

\begin{equation}\label{cc}
 \sum _{p\in P}\sg (p)\prod _{j=1}^{n}\frac{f(\langle \Lambda _{j}(p),
{\bf x}\rangle)}{\langle \Lambda _{j}(p), {\bf x}\rangle}= [M^{2n}]
+ \sum _{\| \omega \| >0}[G_{\omega}(M^{2n})](ch_{U}{\bf u})^{\omega} \; .
\end{equation}

If in the left hand side of this equation we put $t{\bf x}$ instead
of ${\bf x}$ and then multiplying it with $t^{n}$ we obtain the
following result.

\begin{prop}\label{cobc}
The coefficient for $t^{n}$ in the series in $t$
\[
\sum _{p\in P}\sg (p)\prod _{j=1}^{n}\frac{f(t\langle \Lambda
_{j}(p), {\bf x}\rangle)}{\langle \Lambda _{j}(p), {\bf x}\rangle}
\]
represents the complex cobordism class $[M^{2n}]$.
\end{prop}

\begin{prop}\label{coeffzero}
The coefficient for $t^{l}$ in the series in $t$
\[
\sum _{p\in P}\sg (p)\prod _{j=1}^{n}\frac{f(t\langle \Lambda
_{j}(p), {\bf x}\rangle)}{\langle \Lambda _{j}(p), {\bf x}\rangle}
\]
is equal to zero for $0\leqslant l\leqslant n-1$.
\end{prop}

\section{Torus action on homogeneous spaces with positive Euler characteristic. }

Let $G/H$ be a compact homogeneous space of positive Euler
characteristic. It means that $G$ is a compact connected Lie group
and $H$ its connected closed subgroup, such that $\rk G = \rk H$.
Let $T$ be the maximal common torus for $G$ and $H$. There is
canonical action $\theta$ of $T$ on $G/H$ given by $t(gH)=(tg)H$,
where $t\in T$ and $gH\in G/H$. Denote by $N_{G}(T)$ the normalizer
of the torus $T$ in $G$. Then $W_{G} = N_{G}(T)/T$ is the Weyl group
for $G$. For the set  of the fixed points  for the action $\theta$
we prove the following.

\begin{prop}\label{fixed}
The set of fixed points under the canonical action $\theta$ of $T$ on
$G/H$ is given by $(N_{G}(T))~\cdot~H$.
\end{prop}

\begin{proof}
It is easily to see  that $gH$ is fixed point for $\theta$ for any
$g\in N_{G}(T)$. If $gH$ is the fixed point under the canonical
action of $T$ on $G/H$ then $t(gH) = gH$ for all $t\in T$. It
follows that $g^{-1}tg \in H$ for all $t\in T$, i.~e.~ $g^{-1}Tg
\subset H$. This gives that $g^{-1}Tg$ is a maximal torus in $H$
and, since any two maximal toruses in $H$ are conjugate, it follows
that $g^{-1}Tg = h^{-1}Th$ for some $h\in H$. Thus,
$(gh)^{-1}T(gh)=T$ what means that $gh\in N_{G}(T)$. But,
$(gh)H=gH$, what proves the statement.
\end{proof}

Since $T \subset N_{G}(T)$ leaves $H$ fixed, the following Lemma is
direct implication of the  Proposition~\ref{fixed}.

\begin{lem}\label{wfixed}
The set of fixed points under the canonical action $\theta$ of $T$ on
$G/H$ is given by $W_{G}\cdot H$.
\end{lem}

Regarding the number of fixed points, it holds the following.

\begin{lem}\label{number}
The number of fixed points under the canonical action $\theta$ of $T$
on $G/H$ is equal to the Euler characteristic $\chi (G/H)$.
\end{lem}

\begin{proof}
Let $g, g^{'}\in N_{G}(T)$ are representatives of the same fixed
point. Then $g^{'}g^{-1}\in H$ and $g^{-1}Tg = T =
(g^{'})^{-1}Tg^{'}$, what gives that $g^{'}g^{-1}Tg(g^{'})^{-1} = T$
and, thus,  $g^{'}g^{-1}\in N_{H}(T)$. This implies that the number
of fixed points is equal to
\begin{equation}
\Big{\|} \frac{N_{G}(T)}{N_{H}(T)} \Big{\|} = \frac{\| \frac{N_{G}
(T)} {T}\|}{\| \frac{N_{H}(T)}{T} \|} = \frac{\| W_{G}\|}{\| W_{H}
\|}  = \chi (G/H) \; .
\end{equation}

The last equality is classical result related to equal ranks
homogeneous spaces, see~\cite{Onishchik}.
\end{proof}

\begin{rem}
The proof of the Lemma~\ref{number} gives that the set of fixed
points under the canonical action $\theta$  of $T$ on $G/H$ can be
obtained as an orbit of $eH$ by an  action of the Weyl  group
$W_{G}$ up to an action of the Weyl group $W_{H}$.
\end{rem}

\section{The weights at the fixed points.}

Denote  by $\gg$, $\hh$ and $\TT$ the Lie algebras for $G$, $H$ respectively and
$T=T^k$, where $k=\rk G=\rk H$. Let $\alpha
_{1},\ldots, \alpha _{m}$ be the roots for $\gg$ related to $\TT$,
where $\dim G=2m+k$.  Recall that the roots for $\gg$ related to
$\TT$ are the weights for the adjoint representation $\Ad _{T}$  of
$T$ which is given  with $\Ad _{T}(t) = d_{e}\ad (t)$, where $\ad
(t)$ is inner automorphism of $G$ defined by the element $t\in T$.
One can always choose the roots for $G$  such that $\alpha
_{n+1},\ldots, \alpha _{m}$ gives the roots for $\hh$ related to
$\TT$, where $\dim H=2(m-n)+k$. The roots $\alpha _{1},\ldots,
\alpha  _{n}$ are called the {\it complementary} roots for $\gg$
related to $\hh$. Using root decomposition for $\gg$ and $\hh$ it
follows that $T_{e}(G/H) \cong \gg _{\alpha _{1}}^{\C}\oplus \ldots
\oplus \gg _{\alpha _{n}}^{\C}$, where by $\gg _{\alpha _{i}}$ is
denoted the root subspace defined with the root $\alpha _{i}$ and
$T_{e}(G/H)$ is the tangent space for $G/H$ at the $e\cdot H$.
 It is obvious
that $\dim _{\R} G/H=2n$.

\subsection{Description of the invariant almost complex structures.}
Assume we are given  an invariant almost complex structure $J$ on
$G/H$. This means that $J$ is invariant under the canonical action
of $G$ on $G/H$. Then according to the paper~\cite{BH}, we can say
the following.
\begin{itemize}
\item Since $J$ is invariant it commutes with adjoint representation
$Ad_{T}$ of the torus T. This implies that $J$ induces the complex
structure on each complementary root subspace $\gg _{\alpha
_{1}},\ldots ,\gg _{\alpha _{n}}$. Therefore, $J$ can be completely
described by the root system $\varepsilon _{1} \alpha _{1},\ldots
,\varepsilon _{n}\alpha _{n}$, where we take $\varepsilon _{i}= \pm
1$ depending  if $J$ and adjoint representation $Ad_{T}$ define the
same orientation on $\gg _{\alpha _{i}}$ or not, where $1\leqslant
i\leqslant n$. The roots $\varepsilon _{k}\alpha _{k}$ are called
{\it the roots of the almost complex structure} $J$.
\item If we assume $J$ to be integrable, it follows that it can be
chosen an ordering on the canonical coordinates of $\TT$ such that
the roots  $\varepsilon _{1} \alpha _{1},\ldots ,\varepsilon _{n}
\alpha _{n}$ which define $J$ make the closed system of positive
roots.
\end{itemize}

Let us assume that $G/H$ admits an invariant almost complex
structure. Consider an isotropy representation $I_{e}$ of $H$ in
$T_{e}(G/H)$ and let it decomposes into $s$ {\it real irreducible
representations} $I_e = I_e^1 + \ldots +I_e^s$. Then it is proved
in~\cite{BH} that  $G/H$ admits exactly $2^s$ invariant almost
complex structures. Because of completeness we recall the proof of
this fact shortly here. Consider the decomposition of $T_{e}(G/H)$
\[
T_{e}(G/H) = \II _{1}\oplus \ldots \oplus\II _{s} 
\]
such that the restriction of $I_e$ on $\II _{i}$ is $I_{e}^i$. The
subspaces $\II _1,\ldots ,\II_{s}$ are  invariant under $T$ and
therefore each of them is the sum of some root subspaces, i.e. $\II
_{i} = \gg _{\alpha _{i_1}}\oplus\ldots \oplus \gg _{\alpha
_{i_j}}$, for some complementary roots $\alpha _{i_1},\ldots ,\alpha
_{i_j}$. Any linear transformation that commutes with $I_e$ leaves
each of $\II _i$ invariant. Since, by assumption $G/H$ admits an
invariant almost complex structure, we have at least one linear
transformation without real eigenvalue that commutes with $I_e$.
This implies that the commuting field for each of $I_{e}^i$ is the
field of complex numbers and, thus, on each $\II _{i}$ we have
exactly two invariant complex structures.

\begin{rem}
Note that this consideration shows that the numbers $\varepsilon
_1,\ldots ,\varepsilon _n$ that define an invariant almost complex
structure may not vary independently.
\end{rem}

\begin{rem}
In this paper we consider almost complex structures on $G/H$ that
are invariant under the canonical action of the group $G$, what, as
we remarked, imposes some relations on $\varepsilon _1,\ldots
,\varepsilon _n$. If we do not require $G$-invariance, but just
$T$-invariance, we will have more degrees of freedom on $\varepsilon
_1,\ldots ,\varepsilon _n$. This paper is going to have continuation,
where, among the other, the case of $T$-invariant structures will be
studied.
\end{rem}

\begin{ex}
Since the isotropy representation for $\C P^{n}$ is irreducible over
reals, it follows that on $\C P^{n}$ we have only two invariant
almost complex structures, which are actually the standard complex
structure and its conjugate.
\end{ex}

\begin{ex}
The flag manifold $U(n)/T^{n}$ admits $2^{m}$ invariant almost
complex structure, where $m=\frac{n(n-1)}{2}$. By~\cite{BH} only two
of them, conjugate to each other, are integrable.
\end{ex}

\begin{ex}
As we already mentioned, the $10$-dimensional manifold
$M^{10}=U(4)/(U(1)\times U(1)\times U(2))$ is the first example of
homogeneous space, where we have an existence of two non-equivalent
invariant complex structures, see~\cite{BH}. We will in the last
section of this paper also describe cobordism class of $M^{10}$ for
these structures.
\end{ex}

\subsection{The weights at the fixed points.}
We fix now an invariant almost complex structure $J$ on $G/H$ and we
want to describe the weights of the canonical action $\theta$ of $T$
on $G/H$ at the fixed points of this action. If $gH$ is the fixed
point for the action $\theta$, then we obtain a linear map
$d_{g}\theta (t)\colon T_{g}(G/H)\to T_{g}(G/H)$ for all $t\in T$.
Therefore, this action gives rise to the complex representation
$d_{g}\theta$ of $T$ in $(T_{g}(G/H), J)$.

The weights for this representation at identity fixed point  are
described in~\cite{BH}.

\begin{lem}
The weights for the representation $d_{e}\theta$ of $T$ in
$(T_{e}(G/H),J)$ are given by the roots of  an invariant almost
complex structure $J$.
\end{lem}
\begin{proof}
Let us, because of clearness, recall the proof. The inner
automorphism $\ad (t)$, for $t\in T$ induces the map $\overline{\ad
}(t) : G/H\to G/H$ given with $\overline{\ad }(t)(gH) = t(gH)t^{-1} =
(tg)H$. Therefore, $\theta (t) =  \overline{\ad} (t)$ and, thus,
$d_{e}\theta (t) = d_{e}\overline{\ad } (t)$ for any $t\in T$. This
directly gives that the weights  for $d_{e}\theta$ in $(T_{e}(G/H),
J)$ are the roots that define $J$.
\end{proof}

For an arbitrary fixed point we prove the following.

\begin{thm}\label{weights}
Let $gH$ be the fixed point for the canonical action $\theta$ of $T$
on $G/H$. The weights of the induced representation $d_{g}\theta $
of $T$ in $(T_{g}(G/H), J)$ can be obtained from the weights of the
representation $d_{e}\theta$ of $T$ in $(T_{e}(G/H), J)$ by the
action of the Weyl group $W_{G}$ up to the action of the Weyl group
$W_{H}$.
\end{thm}

\begin{proof}
Note  that Lemma~\ref{wfixed} gives that an arbitrary fixed point
can be written as $\rw H$ for some $\rw \in W_{G}/W_{H}$. Fix $\rw
\in W_{G}/W_{H}$ and denote by $l(\rw )$ the action of $\rw$ on
$G/H$, given by $l(\rw)gH=(\rw g)H$ and by $\ad (\rw )$ the inner
automorphism of $G$ given by $\rw$.

We observe that $\theta \circ \ad (\rw )  = \ad (\rw )\circ \theta
$. Then $d_{e}\theta \circ d_{e}\ad (\rw ) = d_{e}\ad (\rw )\circ
d_{e}\theta$. This implies that  the weights for $d_{e}\theta \circ
d_{e}\ad (\rw )$ we get by the action of $d_{e}\ad (\rw)$ on the
weights for $d_{e}\theta$. From the other hand $\theta (\ad (\rw
)t)gH = ({\rw}^{-1}t\rw g)H = (l({\rw}^{-1})\circ \theta (t)\circ
l(\rw ))gH$ what implies that $d_{e}(\theta \circ \ad (\rw )) =
d_{\rw}l({\rw}^{-1})\circ d_{\rw}\theta \circ d_{e}l(\rw )$. This
gives that if, using the map $d_{\rw}l_{{\rw}^{-1}}$, we lift the
weights for $d_{\rw}\theta$ from $T_{\rw}(G/H)$ to $T_{e}(G/H)$, we
get that they coincide with the weights for $d_{e}\theta \circ
d_{e}\ad (\rw)$. Therefore, the weights for $d_{\rw}\theta$ we can
get by the an action of the element $\rw$ on the weights for
$d_{e}\theta$.
\end{proof}

\section{The cobordism classes of homogeneous spaces with positive
Euler characteristic}

\begin{thm}
Let $G/H$ be homogeneous space of compact connected Lie group such
that $\rk G=\rk H = k$ and $\dim G/H=2n$ and consider the canonical
action $\theta$ of maximal torus $T=T^k$ for $G$ and $H$  on $G/H$.
Assume we are given an invariant almost complex structure $J$ on
$G/H$. Let  $\Lambda_{j}=\varepsilon _{j} \alpha _{j}$, $1\leqslant j\leqslant
n$, where $\varepsilon _{1} \alpha _{1},\ldots ,\varepsilon
_{n}\alpha _{n}$ are the complementary roots of $G$ related to $H$
which define an invariant almost complex structure $J$. Then the
toric genus for $(G/H, J)$ is given with
\begin{equation}
\Phi (G/H, J) = \sum _{\rw \in W_{G}/W_{H}} \prod _{j=1}^{n}\frac{1}
{ [\rw (\Lambda_{j})]({\bf u})} \ .
\end{equation}
\end{thm}

\begin{proof}
Rewriting the  Theorem~\ref{EGFP}, since all fixed points have sign
$+1$, we get that the toric genus for $(G/H, J)$  is
\begin{equation}\label{EQH}
\Phi (G/H, J) = \sum _{p\in P} \prod _{j=1}^{n}\frac{1}{[\Lambda
_{j}(p)]({\bf u})} \; ,
\end{equation}
where $P$ is the set of isolated fixed points and $(\Lambda
_{1}(p),\ldots ,\Lambda _{n}(p))$ is the weight vector of the
representation for $T$  in  $T_{p}(G/H)$  associated to an action
$\theta$. By Theorem~\ref{number}, the set of fixed points $P$
coincides with the orbit of the action of $W_{G}/W_{H}$ on $eH$ and
also by Theorem~\ref{weights} the set of weight vectors at fixed
points coincides with the orbit of the action of $W_{G}/W_{H}$ on
the weight vector $\Lambda$ at $eH$.  The result follows if we put
this data into formula~\eqref{EQH}.
\end{proof}

\begin{cor}
The Chern character of the toric genus for homogeneous space $(G/H,
J)$ is given with
\begin{equation}\label{CC}
ch_{U}\Phi(G/H, J) = \sum _{\rw \in W_{G}/W_{H}}\prod _{j=1}^{n}
\frac{f (\langle \rw (\Lambda_{j}), {\bf x}\rangle )}{\langle \rw
(\Lambda_{j}),{\bf x}\rangle} \; ,
\end{equation}
where $f(t) = 1+ \sum\limits_{i\geqslant 1}a_{i}t^{i}$ for $a_{i}
\in \Omega _{U}^{-2i}(\Z), \; {\bf x}=(x_1,\ldots ,x_{k})$ and by
$\langle \Lambda _{j}, {\bf x} \rangle=\sum\limits_{l=1}^k
\Lambda_j^lx_l$ is denoted the weight vector $\Lambda _{j}$ of
$T^k$-representation at $e\cdot H$.
\end{cor}

\begin{cor}\label{chom}
The cobordism class for $(G/H, J)$ is given as the coefficient for
$t^{n}$ in the series in $t$
\begin{equation}
\sum _{\rw \in W_{G}/W_{H}}\prod _{j=1}^{n} \frac{f(t\langle \rw
(\Lambda_{j}), {\bf x}\rangle )}{\langle \rw (\Lambda_{j}) ,{\bf
x}\rangle} \; .
\end{equation}
\end{cor}

\begin{rem}
Since the weights of different invariant almost complex structures
on the fixed homogeneous space differ only by sign, the
Corollary~\ref{chom} provides the way for comparing cobordism
classes of two such structures without having their cobordism
classes explicitly computed.
\end{rem}

\section{Characteristic numbers of homogeneous spaces
with positive Euler characteristic.}

\subsection{Generally about stable complex manifolds.} Let $M^{2n}$ be
tangentially stable complex manifold whose given action $\theta$ of
the torus $T^{k}$ on $M^{2n}$ has only isolated fixed points. Denote
by $P$ the set of fixed points for $\theta$ and set $t_{j}(p) =
\langle \Lambda _{j}(p), {\bf x}\rangle$, where $\{ \Lambda _{j}(p),
\; j=1,\ldots,n\}$ are the weight vectors of the representation of
$T^{k}$ at a fixed point $p$ given by the action $\theta$ and ${\bf
x}=(x_{1},\ldots x_{k})$.

Set
\begin{equation}\label{fdecomp}
\prod_{i=1}^{n}f(t_i)=1+\sum f_\omega(t_1, \ldots, t_n)a^\omega \; .
\end{equation}

Using this notation the Proposition~\ref{coeffzero} could be
formulated in the following way.
\begin{prop}\label{coeffzero1}
 For any $\omega$ with  $0\leqslant \|\omega\| \leqslant (n-1)$ we have that
\begin{equation}
\sum _{p\in P}\sg (p)\cdot \frac{f_\omega(t_1(p), \ldots,
t_n(p))}{t_1(p) \cdots t_n(p)} = 0 \; .
\end{equation}
\end{prop}

Note that the Proposition~\ref{coeffzero1} gives the strong
constraints on the set of signs $\{\sg (p)\}$ and the set of weights
$\{\Lambda _{j}(p)\}$ at fixed points for some tangentially stable
complex manifold.  For example $\|\omega\|=0$ and $\|\omega\|=1$
gives that the signs and the weights at fixed points have to satisfy
the following relations.

\begin{cor}
\begin{equation}
\sum _{p\in P}\sg (p)\cdot \frac{1}{t_1(p) \cdots t_n(p)} = 0 \; ,
\end{equation}
\begin{equation}
\sum _{p\in P}\sg (p)\cdot \frac{\sum\limits _{i=1}^{n}t_{j}(p)}{t_1(p)
\cdots t_n(p)} = 0 \; .
\end{equation}
\end{cor}

As we already mentioned in~\eqref{charact} the cobordism class for
$M^{2n}$ can be represented as

\begin{equation}
[M^{2n}]= \sum_{\| \omega \|= n}s_{\omega}(\tau (M^{2n}))a^{\omega} \; ,
\end{equation}
where $\omega = (i_{1},\ldots ,i_{n})$, $\| \omega \| = \sum _{l=1}^{n}
l\cdot i_{l}$ and $a ^{\omega} = a_{1}^{i_1}\cdots a_{n}^{i_{n}}$.

If the given action $\theta$ of $T^{k}$ on $M^{2n}$ is with isolated
fixed points, the coefficients $s_{\omega}(\tau (M^{2n}))$ can be
explicitly described using Proposition~\ref{cobc} and
expression~\eqref{fdecomp}.

\begin{thm}\label{sthm}
Let $M^{2n}$ be tangentially stable complex manifold whose given
action $\theta$ of the $T^{k}$ have only isolated fixed points.
Denote by $P$  the set of fixed points for $\theta$ and set
$t_{j}(p) = \langle \Lambda _{j}(p), {\bf x}\rangle$, where $\Lambda
_{j}(p)$ are the weight vectors of the representation of $T^{k}$
given by the action $\theta$ and ${\bf x}=(x_{1},\ldots x_{k})$.
Then for $\|\omega\|=n$
\begin{equation}\label{somega}
s_{\omega}(\tau (M^{2n}))=\sum _{p\in P}\sg (p)\cdot
\frac{f_\omega(t_1(p), \ldots, t_n(p))}{t_1(p) \cdots t_n(p)} \; .
\end{equation}
\end{thm}

\begin{ex}\label{ex2}
\[ s_{(n,\ldots,0)}(\tau (M^{2n}))=
\sum _{p\in P}\sg (p)\; . \]
\end{ex}

\begin{ex}\label{ex3}
\[ s_{(0,\ldots,1)}(\tau (M^{2n}))=s_n(M^{2n})=
\sum _{p\in P}\sg (p) \frac{\sum\limits_{j=1}^n t_j^n(p)}{t_1(p)
\cdots t_n(p)} \; . \]
\end{ex}

\begin{rem}\label{int}
Note that the left hand side of~\eqref{somega} in the
Theorem~\ref{sthm}  is an  integer number $s_{\omega}(\tau
(M^{2n}))$ while the right hand side is a rational function in
variables $x_1, \ldots, x_k$. So this theorem imposes strong
restrictions on the sets of signs $\{\sg (p)\}$ and weight vectors
$\{\Lambda_j (p)\}$ of the fixed points.
%See for example the case $\mathbb{C}P^n$ below.
\end{rem}

\subsection{Homogeneous spaces of positive Euler characteristic
and with invariant almost complex structure.} Let us assume
$M^{2n}$ to be homogeneous space $G/H$ of positive Euler
characteristic with canonical action of a maximal torus and endowed
with an invariant almost complex structure $J$. All fixed points
have sign $+1$ and taking into account the Theorem~\ref{weights}, the
Proposition~\ref{coeffzero1} gives that the weights at the fixed
points have to satisfy the following relations.

\begin{cor}
For any $\omega$ with $0\leqslant \|\omega\| \leqslant (n-1)$ where
$2n = \dim G/H$ we have that
\begin{equation}
\sum _{\rw \in W_{G}/W_{H}} \rw \Big( \frac{f_\omega(t_1, \ldots,
t_n)}{t_1 \cdots t_n}\Big ) = 0 \; ,
\end{equation}
where $t_j = \langle \Lambda _{j}, {\bf x}\rangle$ and $\Lambda
_{j}, \; 1\leqslant j\leqslant n$, are the weights  at the fixed
point $e \cdot H$.
\end{cor}

In the same way, the Theorem~\ref{sthm} implies that

\begin{thm}\label{s}
For $M^{2n}=G/H$ and $t_j=\langle  \Lambda_j,{\bf x} \rangle$, where
$\langle  \Lambda_j,{\bf x} \rangle=\sum\limits_{l=1}^k \Lambda_j^l
x_l, \; {\bf x}=(x_1,\ldots,x_k), \; k=\rk G=\rk H$, we have
\begin{equation}
s_{\omega}(\tau (M^{2n}))=\sum _{\rw \in W_{G}/W_{H}}\rw \Big(
\frac{f_\omega(t_1, \ldots, t_n)}{t_1 \cdots t_n}\Big) 
\end{equation}
for any $\omega$ such that $\|\omega\|=n$.
\end{thm}

\begin{ex}\label{euler}
\[
s_{(n,\ldots,0)}(G/H, J)= \| W_{G}/W_{H}\| = \chi (G/H) 
\]
and, therefore, $s_{(n,\ldots ,0)}(G/H, J)$ does not depend on
invariant almost complex structure $J$.
\end{ex}

\begin{cor}\label{sn}
\[
s_{(0,\ldots,1)}(G/H, J)= s_n(G/H, J))= \sum _{\rw \in W_{G}/W_{H}}
\rw \Big( \frac{\sum\limits_{j=1}^n t_j^n}{t_1 \cdots t_n}\Big) \; . \]
\end{cor}

\begin{ex}In the case $\mathbb{C}P^n=G/H$ where $G=U(n+1),\;H=U(1)\times
U(n)$ we have action of $T^{n+1}$ and related to the standard
complex structure the weights are given with $\langle\Lambda_j,{\bf
x}\rangle = x_j-x_{n+1}, \;j=1,\ldots, n$ and
$W_{G}/W_{H}=\mathbb{Z}_{n+1}$ is cyclic group. So

\begin{equation}
s_n(\mathbb{C}P^n) = \sum_{i=1}^{n+1}\frac{\sum\limits_{j\neq i}
(x_i-x_j)^n}{\prod\limits_{j\neq i}(x_i-x_j)}= n+1 \; .
\end{equation}
\end{ex}

\begin{ex}
Let us consider Grassmann manifold $G_{q+2,2}=G/H$ where
$G=U(q+2),\;H=U(2)\times U(q)$. We have here the canonical action of
the torus $T^{q+2}$. The weights for this action at identity point
related to the standard complex structure are given with $\langle
\Lambda _{ij}, {\bf x}\rangle = x_i - x_j$, where $i=1,2$ and
$3\leqslant j\leqslant q+2$. There are $\|W_{U(q+2)}/W_{U(2)\times
U(q)}\|=\frac{(q+2)(q+1)}{2}$ fixed points for this action.
Therefore

\begin{equation}\label{gn2}
s_{2q}(G_{q+2,2}) = \sum _{\rw \in W_{U(q+2)}/W_{U(2)\times U(q)}}
\rw \Big( \frac{\sum\limits_{j=3}^{q+2}\big(
(x_1-x_j)^{2q}+(x_2-x_j)^{2q}\big)
}{\prod\limits_{j=3}^{q+2}(x_1-x_j)(x_2-x_j)}\Big) \; .
\end{equation}

The action of the group  $W_{U(q+2)}/W_{U(2)\times U(q)}$ on the
weights at the identity  point in formula~\eqref{gn2} is given by
the following permutations $\rw=\rw _{kl}$:
\begin{align*}
\rw _{00} &= Id \; ,\\
\rw _{k0}(1) &= k, \; \rw _{k0}(k) =1, \; \mbox{where} \;3\leqslant
k\leqslant q+2 \; ,\\
\rw _{0l}(2) &= l, \; \; \rw _{0l}(l)\; = 2, \; \mbox{where} \;
3\leqslant l\leqslant q+2 \; ,
\end{align*}
\[ \rw _{kl}(1)=k, \;  \rw _{kl}(k) = 1, \; \rw _{kl}(2)=l, \; \rw
_{kl}(l) =2 \;\; \mbox{for}\; 3\leqslant k\leqslant q+1, \;
k+1\leqslant l\leqslant q+2 \; . \]

As we remarked before (see Remark~\ref{int}), the expression on the
right hand side in~\eqref{gn2} is an integer, so we can get a value
for $s_{2q}$ by choosing the appropriate values for the vector
$(x_1,\ldots ,x_{q+2})$. For example, if we take $q=2$ and
$(x_1,x_2,x_3,x_4) = (1,2,3,4)$ the straightforward application of
formula~\eqref{gn2} will give that $s_{4}(G_{4,2}) = -20$.
\end{ex}

\begin{ex}In the case $G_{q+l,l}=G/H$ where $G=U(q+l),\; H=U(q)\times
U(l)$ we have
\begin{equation}
s_{lq}(G_{q+l,l}) = \sum_{\sigma \in S_{q+l}/(S_q \times S_l)}\sigma
\Big( \frac{\sum (x_i-x_j)^{lq}}{\prod(x_i-x_j)}\Big) \; ,
\end{equation}
where $1 \leqslant i \leqslant q, \; (q+1) \leqslant j \leqslant
(q+l)$ and $S_{q+l}$ is the symmetric group.
\end{ex}
We consider later, in the Section~\ref{app}, the case of this
Grassmann manifold in more details.
\subsubsection{Chern numbers.}
We want to deduce an explicit relations between cohomology
characteristic numbers $s_{\omega}$ and classical Chern numbers for
an invariant almost complex structure on $G/H$.

\begin{prop}\label{orbit}
The number $s_{\omega}(\tau (M^{2n}))$, where $\omega=(i_1, \ldots,
i_n), \; \|\omega\|=n$, is the characteristic number that
corresponds to the characteristic class given by the orbit of the
monomial
\[
(u_{1}\cdots u_{i_1})(u_{i_1+1}^2\cdots u_{i_1+i_2}^2)\cdots
(u_{i_1+\ldots i_{n-1}+1}^{n}\cdots u_{i_1+\ldots i_n}^{n}) \; .
\]
\end{prop}

\begin{rem}
Let $\xi=(j_1,\ldots, j_n)$ and ${\bf u}^\xi=u_{1}^{j_1}\cdots
u_{n}^{j_n}$. The orbit of the monomial ${\bf u}^\xi$ is defined
with
\[
O({\bf u}^\xi) = \sum {\bf u}^{\xi'} \; ,
\]
where the sum is over the orbit $\{\xi'=\sigma \xi, \; \sigma \in
S_n \}$ of the vector $\xi \in \mathbb{Z}^n$ under the symmetric
group $S_{n}$ acting by permutation of coordinates of $\xi$.
\end{rem} 

\begin{ex}\label{ex1}
If we take $\omega = (n,\ldots ,0)$ we need to compute the
coefficient for $a_{1}^{n}$ and it is given as an orbit
$O(u_{1}\cdots u_{n})$ and that is elementary symmetric function
$\sigma _{n}$. If we take $\omega =(0,\ldots ,1)$ then we should
compute the coefficient for $a_{n}$ and it is given with
$O(u_{1}^{n})=\sum _{j=1}^{n} u_{j}^{n}$, what is Newton polynomial.
\end{ex}

It is well known fact from the algebra of symmetric functions that
the orbits of monomials give the additive basis for the algebra of
symmetric functions. Therefore, any orbit of monomial can be
expressed through elementary symmetric functions and vice versa. It
gives a relation between the characteristic numbers $s_{\omega}$ in
terms of Chern characteristic numbers $c^\omega=c_1^{i_1}\cdots
c_n^{i_n}$ for an almost complex homogeneous space $(G/H, J)$.

\begin{thm}\label{sc}
Let $\omega =(i_1,\ldots, i_n), \; \|\omega\|=n$, and assume that
the orbit of the monomial
\[
(u_{1}\cdots u_{i_1})(u_{i_1+1}^2\cdots u_{i_1+i_2}^2)\cdots
(u_{i_1+\ldots i_{n-1}+1}^{n}\cdots u_{i_1+\ldots i_n}^{n})
\]
 is expressed through the elementary symmetric function as
\begin{equation}\label{change}
O((u_{1}\cdots u_{i_1})(u_{i_1+1}^2\cdots u_{i_1+i_2}^2)\cdots
(u_{i_1+\ldots i_{n-1}+1}^{n}\cdots u_{i_1+\ldots i_n}^{n}))=
\end{equation}
\[
= \sum _{\| \xi \| = n} \beta_{\xi}\sigma _{1}^{l_1}\cdots \sigma
_{n}^{l_n} 
\]
for some $\beta_{\xi}\in \Z$ and $\| \xi\| = \sum _{j=1}^{n}j\cdot
l_j$, where $\xi=(l_1, \ldots, l_n)$. Then  it holds
\begin{equation}\label{nocohomchern}
s_{\omega}(G/H, J) = \sum _{\rw \in W_{G}/W_{H}}\rw \Big(
\frac{f_\omega(t_1, \ldots, t_n)}{t_1 \dots t_n}\Big)= \sum _{\| \xi
\| = n} \beta_{\xi}c_{1}^{l_1}\cdots c_{n}^{l_n} \; ,
\end{equation}
where $c_{i}$ are the Chern classes for the tangent bundle of
$(G/H,J)$.
\end{thm}

\begin{rem}
Let $p(n)$ denote the number of  partitions of the number $n$. By
varying $\omega$, the equation~\eqref{nocohomchern} gives the system
of $p(n)$ linear equations in Chern numbers whose determinant is,
by~\eqref{change},  non-zero. Therefore, it provides the explicit
formulas for the computation of Chern numbers.
\end{rem}

\begin{rem}
We want to point that relation~\eqref{nocohomchern} in the
Theorem~\ref{sc} together with Theorem~\ref{s} proves that the Chern
numbers for $(G/H, J)$ can be computed without having  any
information on cohomology for $G/H$.
\end{rem}

\begin{ex}\label{ex4}
We provide the  direct application of the Theorem~\ref{sc}
following Example~\ref{ex1}. It is straightforward to see  that
$s_{(n,\ldots ,0)}(G/H) = c_{n}(G/H)$ for any invariant almost
complex structure. This together with Example~\ref{euler} gives that
$c_{n}(G/H) = \chi (G/H)$.
\end{ex}

We want to add that it is given in~\cite{ms} a description of the
numbers $s_{I}$ that correspond to our characteristic numbers
$s_{\omega}$, but the numerations $I$ and $\omega$ are  different.
To the partition $i\in I$ correspond the n-tuple $\omega
=(i_1,\ldots ,i_n)$ such that $i_k$ is equal to the number of
appearances of the number $k$ in the partition $i$.

\section{Some applications.}\label{app}

\subsection{Flag manifolds $U(n)/T^n$.}
We consider invariant complex structure on $U(n)/T^n$.
Recall~\cite{Adams} that the
Weyl group $W_{U(n)}$ is symmetric group and it permutes the
coordinate $x_1,\ldots ,x_n$ on Lie algebra $\TT ^{n}$ for $T^{n}$.
The canonical action of the torus $T^{n}$ on this manifold has $\|
W_{U(n)}\| = \chi (U(n)/T^n) = n!$ fixed points and its weights at
identity point are given by the roots of U(n).

We first consider the case $n=3$ and apply our results to explicitly
compute cobordism class and Chern numbers for $U(3)/T^3$. The roots
for U(3) are $x_1-x_2$, $x_1-x_3$ and $x_2-x_3$. Therefore the
cobordism class for $U(3)/T^3$ is given as the coefficient for $t^3$
in the polynomial
\[
[U(3)/T^3] = \sum _{\sigma \in S_3}\sigma
\Big(\frac{f(t(x_1-x_2))f(t(x_1-x_3))f(t(x_2-x_3))}{(x_1-x_2)
(x_1-x_3)(x_2-x_3)}\Big) \; ,
\]
where $f(t)=1+a_1t+a_2t^2+a_3t^3$, what implies
\[
[U(3)/T^3] = 6(a_1^3 +a_1a_2-a_3) \; .
\]
This gives that the characteristic numbers $s_\omega$ for $U(3)/T^3$
are
\[
s_{(3,0,0)}=6, \quad s_{(1,1,0)} = 6, \quad s_{(0,0,1)}=-6 \; .
\]
By the Theorem~\ref{sc} we have the following relations between the
characteristic numbers $s_\omega$ and the Chern numbers $c^\omega$
\[
c_3 = 6, \; c_1c_2 - 3c_3=6, \; c_1^3-3c_1c_2+3c_3=-6, \; \;
\mbox{what gives}\; \; c_1c_2=24, \; c_1^3=48 \; .
\]

To simplify the notations we take further $\Delta_n =
\prod\limits_{1\leqslant i<j\leqslant n}(x_i-x_j)$.

\begin{thm}
The cobordism class for the flag manifold $U(n)/T^n$ is given as the
coefficient for $t^{\frac{n(n-1)}{2}}$ in the series in $t$
\begin{equation}
\frac{1}{\Delta_n}\sum _{\sigma \in S_{n}} \sg (\sigma) \sigma
\Big(\prod\limits_{1\leqslant i<j\leqslant n}f(t(x_i-x_j)\Big) \; ,
\end{equation}
where $f(t) = 1+\sum\limits_{i\geqslant 1}a_{i}t^{i}$ and $\sg
(\sigma )$ is the sign $\pm 1$ of the permutation $\sigma$.
\end{thm}

\subsubsection{Using of divided difference operators.}
Consider the ring of the symmetric polynomials $\Sym_n \subset
\mathbb{Z}[x_1, \ldots, x_n]$. There is a linear operator (see
\cite{Macdonald-95})

\[ L : \mathbb{Z}[x_1, \ldots, x_n] \longrightarrow \Sym_n  : \;
 L{\bf x}^\xi = \frac{1}{\Delta_n }
\sum _{\sigma \in S_{n}} \sg (\sigma) \sigma ({\bf x}^\xi) \; , \] where
$\xi=(j_1, \ldots, j_n)$ and ${\bf x}^\xi = x_1^{j_1}, \ldots,
x_n^{j_n}$.

It follows from the definition of Schur polynomials $\Sh(x_1,
\ldots, x_n)$ where $\lambda=(\lambda_1 \geqslant \lambda_2
\geqslant \cdots \geqslant \lambda_n \geqslant 0)$ (see
\cite{Macdonald-95}), that
\[ L{\bf x}^{\lambda+\delta} = \Sh(x_1,\ldots, x_n) \; , \]
where $\delta=(n-1,n-2,\ldots,1,0)$ and $L{\bf x}^\delta=1$.
Moreover, the operator $L$ have the following properties:
\begin{itemize}
\item $L{\bf x}^\xi=0$, if $j_1 \geqslant j_2 \geqslant \cdots
\geqslant j_n\geqslant 0$ and $\xi \neq \lambda+\delta$ for some
$\lambda=(\lambda_1 \geqslant \lambda_2 \geqslant \cdots \geqslant
\lambda_n \geqslant 0)$;
\item $L{\bf x}^\xi=\sg (\sigma)L\sigma({\bf x}^\xi)$, where $\xi=(j_1,
\ldots, j_n)$ and $\sigma({\bf x}^\xi)={\bf x}^{\xi'} $, where
$\sigma \in S_n$ and $\xi'=(j'_1 \geqslant j'_2 \geqslant \cdots
\geqslant j'_n\geqslant 0)$;
\item $L$ is a homomorphism of $\Sym_n$-modules.
\end{itemize}

We have
\begin{equation}\label{P}
 \prod\limits_{1\leqslant i<j\leqslant n}f(t(x_i-x_j))= 1+\sum_{|\xi|>0}
P_\xi(a_1, \ldots, a_n, \ldots)t^{|\xi|}{\bf x}^\xi \; ,
\end{equation}
 where
$|\xi|=\sum\limits_{q=1}^n j_q$.

\begin{cor} \label{L}
Set $m=\frac{n(n-1)}{2}$. The cobordism class for the flag manifold
$U(n)/T^n$ is given by the formula
\begin{equation}\label{CL}
[U(n)/T^n] = \sum_{|\xi|=m} P_\xi(a_1, \ldots, a_n, \ldots)L{\bf
x}^\xi \; .
\end{equation}
\end{cor}

\begin{rem}
As we will show in Corollary~\ref{cor8} below, polynomials $P_{\xi}$ in the formula~\eqref{CL} appears to be polynomials only in variables $a_1,\ldots, a_{2n-3}$.
\end{rem}

The characteristic number $s_{m}$ for $U(n)/T^{n}$ is given as
\begin{equation} \label{sm}
s_{m}(U(n)/T^n) = \sum _{1\leqslant i<j\leqslant n} L(x_i-x_j)^{m} \;
.
\end{equation}

\begin{rem}\label{delta}
The first property of the operator $L$ gives that for any $\xi$ such
that $|\xi | =m$, we will have  $L{\bf x}^{\xi}=0$, whenever ${\bf
x}^{\xi}\neq \sigma ({\bf x}^{\delta})$ for every $\sigma \in S_n$.
In other words, in order to have $L{\bf x}^{\xi}\neq 0$, we need
${\bf x}^{\xi}$ to contain $n-1$ variables and with different
degrees.
\end{rem}

Remark~\ref{delta} together  with  (\ref{sm}) and Corollary \ref{L},
implies the following:

\begin{cor}
$s_{1}(U(2)/T^2)=2; \; s_{3}(U(3)/T^3)= -6$\; and
\begin{equation}
s_{m}(U(n)/T^n) = 0 \; ,
\end{equation}
where $m=\frac{n(n-1)}{2}$ and $n>3$.
\end{cor}

We can push up this further. Denote by $(u_1,\ldots,u_m) = \big(
(x_i-x_j),\; i<j \big)$, where $m=\frac{n(n-1)}{2}$. Then for
$\omega = (i_1,\ldots ,i_m),\, \|\omega\|=m$ we have that

\begin{equation}
s_{\omega}(U(n)/T^n) = \sum L \Big( (u_1 \cdots u_{i_1})(u_{i_1+1}^2
\cdots u_{i_1+i_2}^2)\cdots \Big) = \sum L x^\xi \; ,
\end{equation}
where $x^\xi=\sigma(x_1^{n-1}x_2^{n-2} \cdots x_{n-1})$ for some
$\sigma \in S_n$.

Therefore, if $\xi=(j_1, \ldots, j_n)$, then
$\max\limits_{p_1,\ldots,p_s}(j_{p_1}+ \cdots+ j_{p_s}) = s\Big(
n-\frac{s+1}{2} \Big), \; 1\leqslant s \leqslant n$. In particular,
it holds that $\max\limits_{p_1,p_2}(j_{p_1}+j_{p_2}) = 2n-3$.

\begin{cor} \label{cor8}
Let $\omega = (i_1,\ldots ,i_m)$ such that $i_k\neq 0$ for some
$k>2n-3$, then
\[ s_{\omega}(U(n)/T^n) = 0 \; .\]
\end{cor}

If $\omega = (i_1,\ldots ,i_k), \,\|\omega\|=m$, does not satisfy
Corollary~\ref{cor8}, but $i_{k_1},\ldots ,i_{k_l} \neq 0$ for some
$k_1,\ldots ,k_l$ then we have that $k_{p}=2(n-1)-q_{p}$, for
$q_p\geqslant 1$, $1\leqslant p\leqslant l$. In this case we can say the following.

\begin{cor}\label{cor9}
If $n\geqslant 2l$ and $\sum\limits_{p=1}^{l}q_{p} < l(2l-1)$ then
\[ s_{\omega}(U(n)/T^n) = 0\; .\]
\end{cor}

\begin{rem}\label{trans}
From the second property of the operator $L$ we obtain that $L
P(x_1, \ldots, x_n)=0$ , whenever $\sigma (P(x_1, \ldots, x_n)) =
\varepsilon P(x_1, \ldots, x_n)$ for a permutation $\sigma \in S_n$,
where $\varepsilon= \pm 1$ and $\varepsilon \cdot \sg(\sigma)=-1$.
This, in particular, gives that $L(P(x_1,\ldots ,x_n)+\sigma
_{ij}(P(x_1,\ldots ,x_n)))=0$ for any transposition $\sigma _{ij}$
of $x_i$ and $x_j$, where $1\leqslant i<j\leqslant n$.
\end{rem}

Using Remark~\ref{trans} we can compute some more characteristic
numbers of flag manifolds.

\begin{cor}\label{cor9}
Let $n=4q$ or $4q+1$ and $\omega=(i_1, \ldots, i_m), \; \| \omega \|
= m$, where $i_{2l-1}=0$ for $l=1, \ldots, \frac{m}{2}$. Then
$s_{\omega}(U(n)/T^n) = 0$.
\end{cor}

Since $\sigma
_{12}((x_1-x_2)^{2l}\prod\limits_{\underset{(i,j)\neq(1,2)}{1\leqslant
i<j\leqslant n}}f(t(x_i-x_j)))=(x_1-x_2)^{2l}
\prod\limits_{\underset{(i,j)\neq(1,2)}{1\leqslant i<j\leqslant
n}}f(t(x_i-x_j))$ we have, also because of Remark~\ref{trans}, that
\[ L \Big(\prod\limits_{1\leqslant i<j\leqslant
n}f(t(x_i-x_j)\Big) = L \Big(\widetilde f(t(x_1-x_2))
\prod\limits_{\underset{(i,j)\neq(1,2)}{1\leqslant i<j\leqslant
n}}f(t(x_i-x_j))\Big)\; ,
\]
where $\widetilde f(t)=\sum\limits_{l\geqslant 1}a_{2l-1}t^{2l-1}$.
Using this property of $L$ once more we obtain

\begin{thm} \label{thm8}
For $n\geqslant 4$ the cobordism class for the flag manifold
$U(n)/T^n$ is given as the coefficient for $t^{\frac{n(n-1)}{2}}$ in
the series in $t$
\begin{equation} \label{tilde}
L \Big(\widetilde f(t(x_1-x_2))\widetilde f(t(x_{n-1}-x_n))
\prod\limits_{\underset{(i,j)\neq(1,2),(n-1,n)}{1\leqslant
i<j\leqslant n}}f(t(x_i-x_j))\Big)\; .
\end{equation}
\end{thm}

\begin{rem}
The Corollary~\ref{cor9} implies that if $s_{\omega}\neq 0$ for some $\omega =(i_1,\ldots ,i_m)$, then for some $1\leqslant l\leqslant \frac{m}{2}$ it has to be $i_{2l-1}\neq 0$. The Theorem~\ref{thm8} gives stronger results that, for $n\geqslant 4$ in polynomials $P_{\xi}$ in~\eqref{CL} each monom contains the product of at least two elements of the form $a_{2i-1}$.
\end{rem}

The Theorem~\ref{thm8} provide a way for direct computation of the
number $s_{\omega}$, for $\omega =(i_1,\ldots ,i_m)$ such that $\|
\omega \| =2$, where $\| \omega \| = i_1+\ldots +i_m$.  For $n>5$ we
have that $s_{\omega}(U(n)/T^n) = 0$  for such $\omega$.  For $n=4$
and $n=5$ these numbers can be computed very straightforward as the
next example shows.

\begin{ex}
We provide
the computation of the characteristic
number $s_{(1,0,0,0,1,0)}$ for $U(4)/T^4$. From the formula
(\ref{tilde}) we obtain immediately:
\begin{align*}
s_{(1,0,0,0,1,0)}(U(4)/T^4) &= L\Big( (x_1-x_2)(x_3-x_4)^5 +
(x_1-x_2)^5(x_3-x_4) \Big)=\\
&=10L \Big( (x_1-x_2)(x_3-x_4)(x_1^2x_2^2+x_3^2x_4^2)\Big)=\\
&=20L \Big( x_1^3x_2^2(x_3-x_4)+(x_1-x_2)x_3^3x_4^2)\Big)=\\
&=40L \Big( x_1^3x_2^2x_3+x_1x_3^3x_4^2\Big)=80\; .
\end{align*}
\end{ex}

\begin{rem}
We want to emphasize that the formula~\eqref{tilde} gives the
description of the cobordism classes of the flag manifolds in terms
of {\it divided difference operators}. The divided difference
operators are defined with (see \cite{BGG})
\[\partial_{ij}P(x_1,\ldots ,x_n)=\frac{1}{x_i-x_j}\Big(P(x_1,\ldots
,x_n)-\sigma _{ij}P(x_1,\ldots ,x_n)\Big),\] where $i<j$. Put
$\sigma _{i,i+1}=\sigma_i$, $\partial_{i,i+1}=\partial_i, \;
1\leqslant i \leqslant n-1$. We can wright down operator $L$ as the
following composition (see \cite{F, Macdonald-91})
\[
L= (\partial_1 \partial_2 \cdots \partial_{n-1})(\partial_1 \partial_2
\cdots \partial_{n-2}) \cdots (\partial_1 \partial_2) \partial_1 \; .
\]

Denote by $\rw _0$ the permutation $(n,n-1, \ldots, 1)$. Wright down
a permutation $\rw \in S_n$ in the form $\rw = \rw _0 \sigma _{i_1}
\cdots \sigma _{i_p}$ and set $\nabla_{\rw} =
\partial_{i_p}\cdots
\partial_{i_1}$. It is natural to set $\nabla_{\rw _0} = I$  --- identity
operator. The space of operators $\nabla_{\rw}$ is dual to the space
of the Schubert polynomials
$\mathfrak{G}_{\rw}=\mathfrak{G}_{\rw}(x_1,\ldots ,x_n)$, since it
follows from their definition that $\mathfrak{G}_{\rw} =
\nabla_{\rw} {\bf x}^\delta$. Note that $\mathfrak{G}_{\rw _0}= {\bf
x}^\delta$. For the identity permutation $e=(1,2,\ldots, n)$ we have
$e=\rw _0 \cdot \rw _0^{-1}$. So $\nabla_e = L$ and
$\mathfrak{G}_{e} = \nabla_{e} {\bf x}^\delta = 1$.

Schubert polynomials were introduced in \cite{BGG} and in \cite{D}
in context of an arbitrary root systems. The main reference on
algebras of operators $\nabla_{\rw}$ and Schubert polynomials
$\mathfrak{G}_{\rw}$ is \cite{Macdonald-91}.

The description of the cohomology rings of the flag manifolds
$U(n)/T^n$and Grassmann manifolds $G_{n,k}=U(n)/(U(k)\times U(n-k))$
in the terms of Schubert polynomials is given in~\cite{F}.

The description of the {\it complex cobordism ring} of the flag manifolds
$G/T$, for $G$ compact, connected Lie group and $T$ its maximal torus,
 in the terms of the Schubert polynomials calculus is given
in~\cite{Bressler-Evens-90, Bressler-Evens-92}.
\end{rem}

\subsection{Grassmann manifolds.} As a next application we will compute
cobordism class, characteristic numbers $s_{\omega}$ and,
consequently, Chern numbers for invariant complex structure on
Grassmannian $G_{4,2} = U(4)/(U(2)\times U(2)) = SU(4)/S(U(2)\times
U(2))$. Note that, it follows by~\cite{BH} that, up to equivalence,
$G_{4,2}$ has one invariant complex structure $J$. The corresponding
Lie algebra description for $G_{4,2}$ is $A_{3}/(\TT ^{1}\oplus A_1
\oplus A_1)$.

The number of the fixed points under the canonical action of $T^{3}$
on $G_{4,2}$ is, by Theorem~\ref{number}, equal to 6. Let
$x_1,x_2,x_3,x_4$ be canonical coordinates on maximal abelian
algebra for $A_3$. Then $x_1, x_2$ and $x_3, x_4$ represents
canonical coordinates for $A_1\oplus A_1$. The weights of this
action at identity point $(T_{e}(G_{4,2}), J)$are given by the
positive complementary roots $x_1-x_3, x_1-x_4, x_2-x_3, x_2-x_4$
for $A_3$ related to $A_1\oplus A_1$ that define $J$.

The Weyl group $W_{U(4)}$ is symmetric group of permutation on
coordinates $x_1,\ldots ,x_4$ and the Weyl group $W_{U(2)\times
U(2)} = W_{U(2)}\times W_{U(2)}$ is the product of symmetric groups
on coordinates $x_1, x_2$ and $x_3, x_4$ respectively. Let
$\rw_{j}\in  W_{U(4)}/W_{U(2)\times U(2)}$. Corollary~\ref{chom}
gives that the cobordism class $[G_{4,2}]$ is the coefficient for
$t^{4}$ in polynomial

\begin{equation}\label{gr}
\begin{split}
\sum _{j=1}^{6}\rw _{j}\Big( \frac{f(t(x_1-x_3))f(t(x_1-x_4))
f(t(x_2-x_3))
f(t(x_2-x_4))}{(x_1-x_3)(x_1-x_4)(x_2-x_3)(x_2-x_4)}\Big) = \qquad \qquad\\
= \frac{1}{4} L\Big( (x_1-x_2)(x_3-x_4)f(t(x_1-x_3))
f(t(x_1-x_4)) f(t(x_2-x_3)) f(t(x_2-x_4)) \Big)\; ,
\end{split}
\end{equation}
where $f(t) = 1 + a_{1}t + a_{2}t^2 + a_{3}t^3 + a_{4}t^4$.

Expanding formula~\eqref{gr} we get that
\[
[G_{4,2}] = 2(3a_{1}^4+ 12a_{1}^2a_{2} + 7a_{2}^2 + 2a_{1}a_{3}-10a_4) \; .
\]

The characteristic numbers $s_{\omega}$ can be read off form this formula:
\[
s_{(4,0,0,0)}=6, \quad s_{(2,1,0,0)}=24, \quad s_{(0,2,0,0)}=14,
\quad s_{(1,0,1,0)}=4, \quad s_{(0,0,0,1)} = -20 \; .
\]

The coefficients $\beta _\xi$ from the Theorem~\ref{sc} can be
explicitly computed and for $8$-dimensional manifold give the
following relation between characteristic numbers $s_{\omega}$ and
Chern numbers:
\[
s_{(0,0,0,1)}= c_{1}^4 - 4c_{1}^2c_{2} + 2c_{2}^2 + 4c_{1}c_{3} - 4c_{4}
, \quad s_{(2,1,0,0)} = c_1c_3 - 4c_4 \; ,
\]
\[
s_{(0,2,0,0)} = c_{2}^2 - 2c_1c_3 + 2c_4, \quad
s_{(1,0,1,0)}=c_1^2c_2-c_1c_3+4c_4-2c_2^2,\quad s_{(4,0,0,0)} = c_4 \; .
\]
We deduce that the Chern numbers for $(G_{4,2}, J)$ are
\[
c_4 = 6, \quad c_1c_3 = 48, \quad  c_2^2 = 98, \quad c_1^2c_2 = 224,
\quad c_1^4 = 512 \; .
\]

The given example  generalizes as follows. Denote by  $\Delta_{p,q}
= \prod\limits_{p\leqslant i<j\leqslant q}(x_i-x_j)$, then
$\Delta_n = \Delta_{1,n}$.
\begin{thm}
The cobordism class for Grassmann manifold $G_{q+l,l}$ is given as
the coefficient for $t^{lq}$ in the series in $t$
\begin{equation}
\sum _{\sigma\in S_{q+l}/S_{q}\times S_{l}}\sigma\Big(\prod
\frac{f(t(x_i-x_j))}{(x_i-x_j)}\Big) = \frac{1}{q!l!}L \Big(
\Delta_q\Delta_{q+1,q+l}\prod f(t(x_i-x_j)) \Big)\; ,
\end{equation}
where $1 \leqslant i \leqslant q, \; (q+1) \leqslant j \leqslant
(q+l)$ and $S_{q+l}$ is the symmetric group.
\end{thm}

\subsection{Homogeneous space $SU(4)/S(U(1)\times U(1)\times U(2))$}

Following~\cite{BH} and~\cite{KT} we know that 10-dimensional space
$M^{10} = SU(4)/S(U(1)\times U(1)\times U(2))$ admits, up to
equivalence, two invariant complex structure $J_1$ and $J_2$ and one
non-integrable invariant almost complex structure $J_{3}$. We
provide here the description of cobordism classes for all of three
invariant almost complex structures. The Chern numbers for all the
invariant almost complex structures are known and they have been
completely computed in~\cite{KT} through multiplication in
cohomology. We provide also their computation
using our method.

The corresponding Lie algebra description for $M^{10}$ is
$A_{3}/(\TT ^{2}\oplus A_1)$. Let $x_1, x_2, x_3, x_4$ be canonical
coordinates on maximal Abelian subalgebra for $A_3$. Then $x_1, x_2$
represent canonical coordinates for $A_1$. The number of fixed
points under the canonical action of $T^{3}$ on $M^{10}$ is, by
Theorem~\ref{number}, equal to 12.

\subsubsection{The invariant complex structure $J_1$.} The weights of the
action of $T^{3}$ on $M^{10}$ at identity point related to $J_1$ are
given by the  complementary roots $x_1-x_3, x_1-x_4, x_2-x_3,
x_2-x_4, x_3-x_4$ for $A_{3}$ related to $A_1$,
(see~\cite{BH},~\cite{KT}). The cobordism class $[M^{1}, J_1]$ is,
by Corollary~\ref{chom}, given as the coefficient for $t^{5}$ in
polynomial

\[
\sum _{j=1}^{12}\rw _{j}\Big( \frac{f(t(x_1-x_3))f(t(x_1-x_4))
f(t(x_2-x_3))f(t(x_2-x_4))
f(t(x_3-x_4))}{(x_1-x_3)(x_1-x_4)(x_2-x_3)(x_2-x_4)(x_3-x_4)}\Big ) \; ,
\]
where $\rw_{j}\in  W_{U(4)}/W_{U(2)}$ and $f(t) = 1 + a_{1}t +
a_{2}t^2 + a_{3}t^3 + a_{4}t^4 + a_{5}t^5$.

Therefore we get that
\[
[M^{10}, J_{1}]= 4(3a_1^5 + 12a_1^3a_2 + 7a_1a_2^2 -
5a_1^2a_3 - 2a_2a_3 - 10a_1a_4 + 5a_5) \; .
\]

Then Theorem~\ref{sc} gives the following relations between
characteristic numbers $s_\omega$ and Chern numbers  for $(M^{10},
J_1)$.

\[
s_{(0,0,0,0,1)} = 20 = c_1^5 - 5c_1^3c_2 + 5c_1^2c_3 + 5c_1c_2^2 -
5c_1c_4 - 5c_2c_3 + 5c_5 \; ,
\]
\[
s_{(1,2,0,0,0)}= 28 = c_2c_3 - 3c_1c_4 + 5c_5, \quad
s_{(2,0,1,0,0)} = -20 = c_1^2c_3 - c_1c_4 - 2c_2c_3 + 5c_5 \; ,
\]
\[
s_{(0,1,1,0,0)} = -8 = -2c_1^2c_3 + c_1c_2^2 - c_2c_3 + 5c_1c_4 - 5
c_5, \quad s_{(3,1,0,0,0)} = 48 = c_1c_4 - 5c_5 \; ,
\]
\[
s_{(1,0,0,1,0)} = -40 = c_1^3c_2 - c_1^2c_3 - 3c_1c_2^2 + c_1c_4 +
5c_2c_3 - 5c_5, \quad s_{(5,0,0,0,0)} = 12 = c_5 \; .
\]

This implies that the Chern numbers for $(M^{10}, J_1)$ are as
follows:
\[
c_5 = 12, \quad c_1c_4 = 108, \quad c_2c_3 = 292, \quad c_1^2c_3 =
612, \quad c_1c_2^2 = 1028, \quad c_1^3c_2 = 2148, \quad c_1^5 =
4500 \; .
\]

\subsubsection{The invariant complex structure $J_2$.} The weights of the
action of $T^{3}$ on $M^{10}$ at identity point related to $J_2$ are
given by the positive complementary roots $x_4-x_1, x_4-x_2,
x_4-x_3, x_1-x_3, x_2-x_3$ for $A_{3}$ related to $A_1$,
(see~\cite{BH},~\cite{KT}). The cobordism class $[M^{1}, J_2]$ is,
by Corollary~\ref{chom}, given as the coefficient for $t^{5}$ in
polynomial

\[
\sum _{j=1}^{12}\rw _{j}\Big( \frac{f(t(x_4-x_1))f(t(x_4-x_2))
f(t(x_4-x_3))f(t(x_1-x_3))
f(t(x_2-x_3))}{(x_4-x_1)(x_4-x_2)(x_4-x_3)(x_1-x_3)(x_2-x_3)}\Big) \; ,
\]
where $\rw_{j}\in  W_{U(4)}/W_{U(2)}$ and $f(t) = 1 + a_{1}t +
a_{2}t^2 + a_{3}t^3 + a_{4}t^4 + a_{5}t^5$.

Therefore we get that
\[
[M^{10}, J_{2}]= 4(3a_1^5 + 12a_1^3a_2 + 7a_1a_2^2 - 5a_1^2a_3 +
8a_2a_3 - 10a_1a_4 - 5a_5) \; .
\]

Applying the same procedure as for above we get that the Chern
numbers for $(M^{10}, J_{2})$ are:
\[
c_5 = 12, \quad c_1c_4 = 108, \quad c_2c_3 = 292, \quad c_1^2c_3 = 612,
\quad c_1c_2^2 = 1068, \quad c_1^3c_2 = 2268, \quad c_1^5 = 4860 \; .
\]

\subsubsection{The invariant almost complex structure $J_3$.} The weights for
the action of $T^{3}$ on $M^{10}$ at identity point related to $J_3$
are given by complementary roots $x_1-x_3, x_2-x_3, x_4-x_1,
x_4-x_2, x_3-x_4$, (see~\cite{KT}). Using Corollary~\ref{chom} we get that the
cobordism class for $(M^{10}, J_3)$ is
\[
[M^{10}, J_{3}] = 4(3a_1^5 - 12a_1^3a_2 + 7a_1a_2^2 + 15a_1^2a_3 -
12a_2a_3 -10a_1a_4 + 15a_5) \; .
\]

The characteristic numbers for $(M^{10}, J_{3})$ are given as
coefficients in its cobordism class, what, as above,  together with
Theorem~\ref{sc} gives that the Chern numbers for $(M^{10}, J_{3})$
are as follows:
\[
c_5 = 12, \quad c_1c_4 = 12, \quad c_2c_3 = 4, \quad c_1^2c_3 = 20,
\quad c_1c_2^2 = -4, \quad c_1^3c_2 = -4, \quad c_1^5 = -20 \; .
\]

\begin{rem}
Further work on the studying of Chern numbers
and the geometry for the generalizations of this example is done
in~\cite{H05} and in~\cite{KT}.
\end{rem}

\bibliographystyle{amsplain}

\end{document}